# A Dual Approach to Triangle Sequences: A Multidimensional Continued Fraction Algorithm

Sami Assaf
University of Notre Dame

Li-Chung Chen
Harvard University

Tegan Cheslack-Postava
Williams College

Benjamin Cooper
Williams College

Alexander Diesl
Johns Hopkins University

Thomas Garrity
Department of Mathematics
Williams College
Williamstown, MA 01267
email:tgarrity@williams.edu *

Mathew Lepinski
Rose-Hulman Institute of Technology

Adam Schuyler
Williams College

**Abstract**

A dual approach to defining the triangle sequence (a type of multidimensional continued fraction algorithm, initially developed in [9]) for a pair of real numbers is presented, providing a new, clean geometric interpretation of the triangle sequence. We give a new criterion for when a triangle sequence uniquely describes a pair of numbers and give the first explicit examples of triangle sequences that do not uniquely describe a pair of reals. Finally, this dual approach yields that the triangle sequence is topologically strongly mixing, meaning in particular that it is topologically ergodic.

*This work was partially supported by the NSF's REU grant to the Williams College SMALL program.

# 1 Introduction

In 1848, Hermite asked Jacobi for methods of expressing a real number as a sequence of integers such that the algebraic properties of the real number are reflected in the periodicity of its sequence. In other words, Hermite wanted a generalization to cubic and higher degree algebraic numbers of the fact that the decimal expansion of a real number is periodic if and only if the real is rational and, more importantly, of the fact that the continued fraction expansion of a real number is periodic if and only if the real is a quadratic irrational. Such attempts are called *multidimensional continued fractions*.

For a good survey of work on multidimensional continued fractions, see Schweiger's *Multidimensional Continued Fractions* [29] (his earlier works [28] and [26] should also be consulted). For many of the algorithms that existed as of 1980, see Brentjes' *Multi-Dimensional Continued Fraction Algorithms* [2]. There is also the deep work of Minkowski [22] [23]. Other work is in [8], [10], [12], [15], [16], [17], [18], [19]. We will be concentrating on understanding the multidimensional continued fraction given in [9].

None of these techniques provides a link that will precisely identify periodicity of integer sequences with algebraic numbers. Almost all of these methods have the property that the periodicity of the sequence will imply algebraticity; none proves the converse. Probably there is no single such technique that will fully answer Hermite's initial question. It is more likely that there is a whole family of techniques, each providing a periodic sequence for different classes of algebraic numbers. For now, each existing method has its own strengths and weaknesses. One way of measuring a method's strength is how many of the properties of traditional continued fractions are generalized by the method.

The method presented in [9] generalizes the Gauss map of the unit interval to a map (called the *triangle map*) of a simplex to itself. This paper shows that the geometric approach to continued fractions (which approximates a line in the plane by better and better integer lattice points) also has an extremely natural interpretation in terms of the triangle map. Further, this method provides a clean description for one of the more interesting features of the triangle map, namely that a given integer sequence need not uniquely describe a point. (This is in marked contrast to most other multidimensional continued fractions.) With this paper's approach, we have a clean description for when a triangle sequence corresponds to a unique point. This description also allows us to determine dynamical properties of the triangle sequence. (Most of the other multidimensional continued fraction algorithms can be shown to have ergodic properties; what prevents us from applying these techniques to the triangle map is the problem of uniqueness, forcing us to develop other techniques.)

We review the relevant facts of continued fractions in section two and of triangle sequences in section three. Section four (which is the start of what is new in this paper) gives a clean description of the vertices of the defining triangles for a given triangle sequence. This leads us in section five to see how the triangle map has a good geometric description in terms of how certain planes move in space about a given ray, in direct analogue to how continued fractions can be defined via adding vectors to get as close as possible to a given ray without crossing the ray. Section six is the longest and most difficult of this paper. The goal of this section is to give a sharp description of precisely when a triangle sequence corresponds to a unique pair of numbers $(\alpha, \beta)$. By section 6.6, enough structure has been developed so that explicit examples of both uniqueness and non-uniqueness can be given. We view the fact that there exists any structure at all as interesting. In section seven, using the machinery developed in the previous section, we show that the triangle map is topologically strongly mixing, which implies, for example, that it is topologically ergodic.

We have developed a Mathematica package for calculating triangle sequences that is available at the web site:

http://www.williams.edu/Mathematics/tgarrity/triangle.html

We would like to thank Lori Pedersen for providing all but the first of the diagrams and for providing many comments. Also, T. Garrity would like to thank the mathematics department at the University of Michigan, where part of this paper was written while he was on sabbatical.

# 2  Continued Fractions

The quickest method for defining the continued fraction expansion for a real number $\alpha \in (0, 1]$ is to use the Gauss map. Set

$$
\begin{aligned}
I &= \{x \mid 0 < x \le 1\} \\
I_k &= \{x \in I \mid k \le \frac{1}{x} < k + 1\}.
\end{aligned}
$$

The *Gauss map* $G : I \to I \cup \{0\}$ is:

$$
G(x) = \frac{1}{x} - k
$$

for $x \in I_k$. Then the continued fraction expansion for any $\alpha \in I$ is the sequence of positive integers $(a_1, a_2, \ldots)$ such that for each $k \ge 0$,

$$
G^k(\alpha) \in I_{a_{k+1}},
$$

where it is understood that if, for some $k$, we have $G^k(\alpha) = 0$, then the continued fraction expansion sequence stops. It is this approach that was directly generalized in [9], where we replaced the unit interval, and its partitioning into subintervals, by a triangle, and its partitioning into subtriangles.

There is a more geometric approach to continued fractions, as explained in [30] on page 187. It is this approach that we will generalize, though as with continued fractions, this approach will yield the same sequence as that in [9]. Given a real number $\alpha \in I$, consider the line $L$ defined by $y = \alpha x$. Define vectors

$$
V_0 = (1, 0) \text{ and } V_{-1} = (0, 1).
$$

Note that these two vectors lie on opposite sides of the line $L$. Define $a_1$ to be the unique positive integer such that the vector

$$V_1 = V_{-1} + a_1 V_0$$

either lies on the line $L$ or on the same side of $L$ as does the vector $V_{-1}$ and the vector

$$V_{-1} + (a_1 + 1)V_0$$

lies on the other side of $L$. If we have constructed vectors $V_{-1}, V_0, V_1, \ldots, V_{n-1}$ such that the odd vectors $V_{2k+1}$ lie on one side of $L$ and the even vectors $V_{2k}$ lie on the other side, then define $a_n$ to be the unique largest positive integer such that the vector

$$V_n = V_{n-2} + a_n V_{n-1}$$

lies on $L$ or on the same side of $L$ as does $V_{n-2}$ but that

$$V_{n-2} + (a_n + 1)V_{n-1}$$

lies on the other side of $L$. If any vector $V_n$ lands on $L$, stop. As shown in [30], this sequence of positive integers $(a_1, a_2, \ldots)$ is the continued fraction expansion of the number $\alpha$. (Note that we do not start with an $a_0$ term, unlike Stark in [30], since we make the initial assumption that the number $\alpha$ is between zero and one.)

## 3   Triangle Sequences

Recall the triangle sequence as developed in [9]. Consider pairs of real numbers $(\alpha, \beta)$ in the triangle $\triangle = \{(x, y) : 1 \geq x \geq y > 0\}$. Partition $\triangle$ into

disjoint triangles

$$\triangle_k = \{(x,y) \in \triangle : 1 - x - ky \geq 0 > 1 - x - (k+1)y\},$$

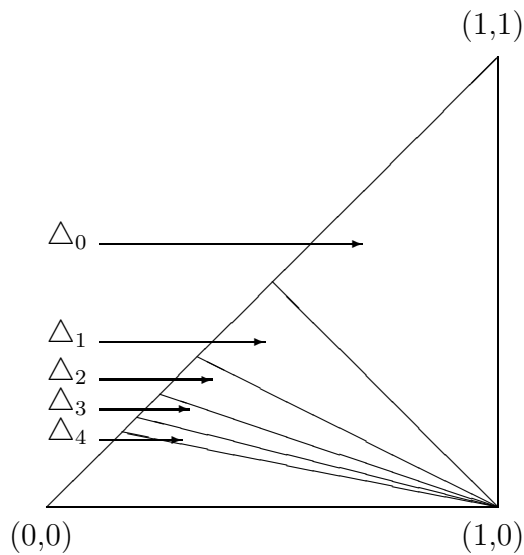

and define $T : \triangle \to \triangle \cup \{(x,0) : 0 \leq x \leq 1\}$ by

$$T(\alpha, \beta) = (\frac{\beta}{\alpha}, \frac{1 - \alpha - k\beta}{\alpha}),$$

if the pair $(\alpha, \beta) \in \triangle_k$. Then the triangle sequence for a pair $(\alpha, \beta)$ will be the infinite sequence of nonnegative integers $(a_0, a_1, a_2, \ldots)$ if $T^k(\alpha, \beta) \in \triangle_{a_k}$. Note that the triangle sequence is said to terminate at step $k$ if $T^k(\alpha, \beta)$ lands on the interval $\{(t, 0) : 0 \leq t \leq 1\}$. As discussed in [9], the hope is that interesting properties of this sequence reflect interesting properties of the original pair $(\alpha, \beta)$. For example, if the sequence is eventually periodic, then both $\alpha$ and $\beta$ are contained in the same cubic number field.

Another way of thinking about triangle sequences is as a method for producing integer lattice vectors in space that approximate the plane $x +$ $\alpha y + \beta z = 0$. Since the normal to this plane is the vector $(1, \alpha, \beta)$, we need to produce vectors whose dot products with $(1, \alpha, \beta)$ are small. We do this inductively as follows. Set

$$C_{-3} = \begin{pmatrix} 1 \\ 0 \\ 0 \end{pmatrix}, C_{-2} = \begin{pmatrix} 0 \\ 1 \\ 0 \end{pmatrix}, C_{-1} = \begin{pmatrix} 0 \\ 0 \\ 1 \end{pmatrix}.$$

If the triangle sequence for $(\alpha, \beta)$ is $(a_0, a_1, a_2, \ldots)$, set

$$C_k = C_{k-3} - C_{k-2} - a_k C_{k-1}.$$

The triangle sequence can in fact be defined in terms of the dot products

$$d_k = (1, \alpha, \beta) \cdot C_k.$$

Assuming we know the number $a_0, \ldots, a_k$, then $a_{k+1}$ is the nonnegative integer such that

$$d_{k-2} - d_{k-1} - a_{k+1} d_k \geq 0 > d_{k-2} - d_{k-1} - (a_{k+1} + 1) d_k.$$

Then

$$d_{k+1} = d_{k-2} - d_{k-1} - a_{k+1} d_k.$$

## 4 Vertices of Triangles

Let $(a_0, a_1, a_2, \ldots)$ be a sequence of nonnegative integers. Define

$$\triangle(a_0, \ldots, a_n) = \{(x, y) : T^k(x, y) \in \triangle(a_k), \text{for all } k \leq n\}.$$

Thus $\triangle(a_0, \ldots, a_n)$ consists of all those points whose first $n + 1$ terms in their triangle sequence are $(a_0, \ldots, a_n)$. As shown in [9], each $\triangle(a_0, \ldots, a_n)$ is indeed a triangle. This section will find a clean formula for the vertices of each of these triangles in terms of the approximating vectors $C_k$.

Define

$$X_k = C_k \times C_{k+1}.$$

By using the recursion formula for the vectors $C_k$, we have by direct calculation

**Proposition 1**

$$X_k = X_{k-3} + a_k X_{k-2} + X_{k-1}$$

Denote each $X_k$ as

$$X_k = \begin{pmatrix} x_k \\ y_k \\ z_k \end{pmatrix}.$$

Then we have by the above formula:

**Corollary 2** *The sequence $\{x_k\}$ is a strictly increasing sequence of positive reals, for $k \geq 0$.*

We need one more piece of notation before we can find the vertices of the triangles $\triangle(a_0, \ldots, a_n)$. For any vectors

$$T = \begin{pmatrix} a \\ b \\ c \end{pmatrix} \text{ and } S = \begin{pmatrix} d \\ e \\ f \end{pmatrix}$$

with $a, d, a + d \neq 0$, define

$$\hat{T} = \begin{pmatrix} \frac{b}{a} \\ \frac{c}{a} \end{pmatrix}$$

and further, define

$$T\hat{+}S = \begin{pmatrix} \frac{b+e}{a+d} \\ \frac{c+f}{a+d} \end{pmatrix}.$$

(Such a sum is called a *Farey sum.*) We can now cleanly describe the vertices for the triangle $\triangle(a_0, \ldots, a_n)$.

**Theorem 3** *The vertices for the triangle* $\triangle(a_0, \ldots, a_n)$ *are* $\hat{X}_{n-1}$, $\hat{X}_n$ *and* $X_n\hat{+}X_{n-2}$.

**Proof:** We do this by induction. The base case is a straightforward calculation. Thus suppose that the vertices for $\triangle(a_0, \ldots, a_{n-1})$ are $\hat{X}_{n-2}$, $\hat{X}_{n-1}$ and $X_{n-1}\hat{+}X_{n-3}$.

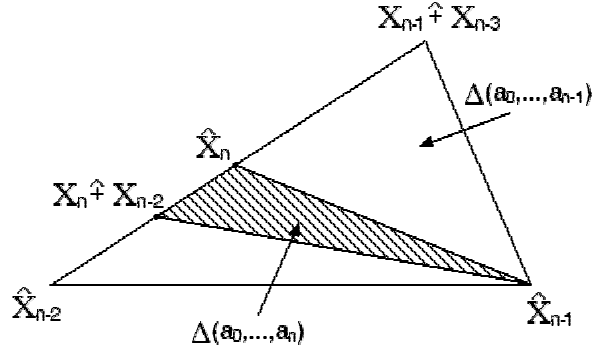

Every point in this triangle has $(a_0, \ldots, a_{n-1})$ as its first $n$ terms in its triangle sequence. Those points $(x, y)$ whose next term will be $a_n$ are those such that

$$(1, x, y) \cdot (C_{n-3} - C_{n-2} - a_n C_{n-1}) \geq 0 > (1, x, y) \cdot (C_{n-3} - C_{n-2} - (a_n + 1)C_{n-1}).$$

We must put this in terms of our conjectured vertices.

Geometrically, in three space with coordinates labeled by $(z, x, y)$, our triangles can be viewed as living in the plane $(z = 1)$. For any vector $T =$ $\begin{pmatrix} a \\ b \\ c \end{pmatrix}$, then the vector $\hat{T}$ can be viewed as the intersection of the ray spanned by $T$ with the plane $(z = 1)$. We can see that, on the segment connecting $\hat{X}_{n-2}$ and $X_{n-1}\hat{+}X_{n-3}$, lie the two points $\hat{X}_n$ and $X_n\hat{+}X_{n-2}$. Let $P$ denote the plane spanned by the vectors $X_{n-1}$ and by $X_n$, and let $Q$ denote the plane spanned by the vectors $X_{n-1}$ and $X_n + X_{n-2}$. In terms of the above diagram, the line segment from $\hat{X}_{n-1}$ to $\hat{X}_n$ is precisely the intersection of the plane $P$ with the triangle $\triangle$, (which, again, is assumed here to be in the plane $(z = 1)$. Likewise, the line segment from $\hat{X}_{n-1}$ to $X_n\hat{+}X_{n-2}$ is the intersection of the plane $Q$ with $\triangle$.

In the first octant, we want to show that the rays spanned by vectors $(1, x, y)$ for points $(x, y) \in \triangle(a_0, \ldots, a_n)$ lie between the planes $P$ and $Q$. By taking cross products of the defining vectors for each plane, note that the normal vectors to the planes $P$ and $Q$ are $C_n = (C_{n-3} - C_{n-2} - a_n C_{n-1})$ and $C_n - C_{n-1} = (C_{n-3} - C_{n-2} - (a_n+1)C_{n-1})$, respectively. Since the basis $X_{n-1}$, $X_n + X_{n-2}$ and $C_n$ has the same orientation as the basis $X_{n-1}$, $X_n + X_{n-2}$ and $C_n - C_{n-1}$, the condition that the ray $(1, x, y)$ is between the planes $P$ and $Q$ in the first octant is precisely that

$$(1, x, y) \cdot (C_{k-3} - C_{k-2} - a_k C_{k-1}) \geq 0 > (1, x, y) \cdot (C_{k-3} - C_{k-2} - (a_k+1)C_{k-1}),$$

which is what we need. $\square$

# 5   The Dual Approach to Triangle Sequences

Given our point $(\alpha, \beta) \in \triangle$, we have constructed a nested sequence of triangles

$$\triangle \supset \triangle(a_0) \supset \triangle(a_0, a_1) \supset \ldots \supset \triangle(a_0, a_1, \ldots, a_n) \supset \ldots.$$

We will see in the rest of this paper that this nested sequence either converges to the initial point $(\alpha, \beta)$ or to a line segment containing $(\alpha, \beta)$.

We now want to see how this provides another clean generalization of continued fractions. Let $(a_0, a_1, \ldots)$ be the continued fraction expansion for a positive real number $\alpha$ and denote the partial convergents by $p_k/q_k = (a_0, \ldots, a_k)$. We have that

$$\begin{pmatrix} q_{k+1} \\ p_{k+1} \end{pmatrix} = \begin{pmatrix} a_{k+1}q_k \\ a_{k+1}p_k \end{pmatrix} + \begin{pmatrix} q_{k-1} \\ p_{k-1} \end{pmatrix} = \begin{pmatrix} a_{k+1}q_k + q_{k-1} \\ a_{k+1}p_k + p_{k-1} \end{pmatrix}.$$

Geometrically, we have, for $v_k = \begin{pmatrix} q_k \\ p_k \end{pmatrix}$

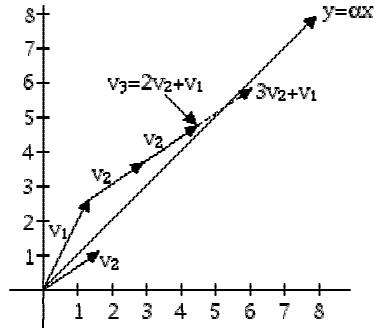

Thus the vectors $\begin{pmatrix} q_{k-1} \\ p_{k-1} \end{pmatrix}$ and $\begin{pmatrix} q_k \\ p_k \end{pmatrix}$ lie on opposite sides of the ray $y = \alpha x$ and $a_{k+1}$ is that positive integer such that the vector

$$a_{k+1} \begin{pmatrix} q_k \\ p_k \end{pmatrix} + \begin{pmatrix} q_{k-1} \\ p_{k-1} \end{pmatrix}$$

lies on the same side of $y = \alpha x$ as $\begin{pmatrix} q_{k-1} \\ p_{k-1} \end{pmatrix}$. These partial convergents also produce for us a nested sequence of intervals $I_1 \supset I_2 \supset \dots$ about the point $\alpha$, where

$$I_{2k} = [p_{2k}/q_{2k}, p_{2k-1}/q_{2k-1}] \text{ and } I_{2k+1} = [p_{2k}/q_{2k}, p_{2k+1}/q_{2k+1}]$$

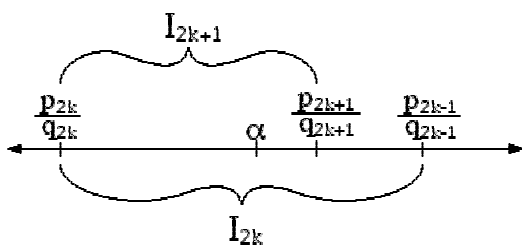

Now to see how our nested sequence of triangles generalizes this. As in the previous section, we put our triangle into the plane $z = 1$. The analogue of the ray $y = \alpha x$ will be the line $(x = \alpha z, y = \beta z)$. Consider the cone through the origin over each triangle $\triangle(a_0, a_1, \dots, a_n)$.

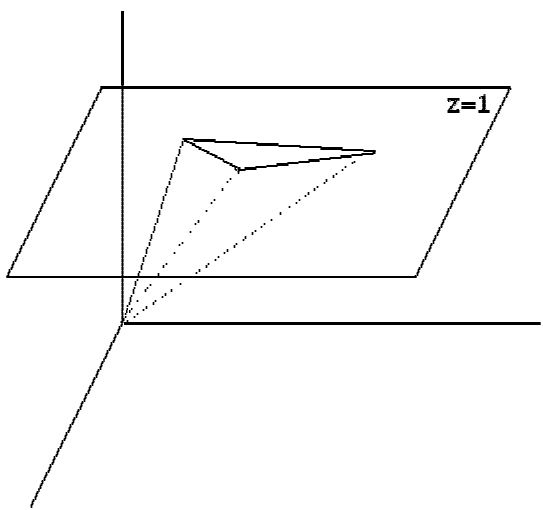

Then the triangle sequence is producing a nested sequence of such cones over the nested sequence of triangles. The analogue of the adding of vectors to get $a_{k+1} \begin{pmatrix} q_k \\ p_k \end{pmatrix} + \begin{pmatrix} q_{k-1} \\ p_{k-1} \end{pmatrix}$ will be planes spanned by vectors $X_{n-1}$ and by the vector $X_{n-1} + X_{n-3} + a_n X_{n-2}$.

Fix a positive integer $n$. For each nonnegative integer $k$, let $P_k$ denote the plane spanned by the vector $X_{n-1}$ and by the vector $X_{n-1} + X_{n-3} + kX_{n-2}$. In the notation from the above proof, we have $P_{a_n} = P$ and $P_{a_n+1} = Q$. Assume we have our pair of numbers $(\alpha, \beta)$ and that we have already found the first $n$ terms of the pair's triangle sequence, $(a_0, a_1, \ldots, a_{n-1})$. We want to see how to find the next term of the triangle sequence solely in terms of the vectors $X_{n-1}$, $X_{n-2}$ and $X_{n-1} + X_{n-3}$. The planes $P_k$ form a family of planes rotating about the ray spanned by $X_{n-1}$, moving away from the plane $P_0$ towards the plane spanned by the vectors $X_{n-1}$ and $X_{n-2}$. Choose $a_n$ to be that positive integer such that the vector $(1, \alpha, \beta)$ lies between the planes $P_{a_n}$ and $P_{a_n+1}$. This is in direct analogue to the geometric development of continued fractions as given in [30].

## 6    Problems with Uniqueness

Triangle sequences have one peculiarity not shared with other multidimensional continued fraction algorithms. Namely, a sequence of non-negative integers need not correspond to a unique pair of real numbers $(\alpha, \beta) \in \triangle$ but could correspond to an entire line segment. The goal of this section is a clean description in terms of the growth of the numbers $a_k$ for when the sequence does correspond to a unique pair $(\alpha, \beta)$. Crudely, we will see that if the terms in the triangle sequence grow sufficiently fast, then we will have non-uniqueness. The existence of a clean criterion for uniqueness and non-uniqueness indicates that the triangle iteration has interesting hidden structure. As an added benefit, the machinery developed here will be critical for our results in section seven on the topological dynamics of the triangle map. As a word of warning, this section is long and detailed.

## 6.1 Parity Results

Before we can address concerns of uniqueness, we need to examine more closely the triangles $\triangle(a_0, \ldots, a_n)$. As shown in [9] (this can also be directly calculated), if two pairs of real numbers are both in some $\triangle(a_0, \ldots, a_n)$, then every point on the line segment connecting the pairs must be in $\triangle(a_0, \ldots, a_n)$. Since the determinant of the Jacobian of the each map $T_k$ is greater than one, this means that only single isolated points or line segments can have the same triangle sequences.

We know that the vertices of the triangle $\triangle(a_0, \ldots, a_n)$ are $\hat{X}_{n-1}$, $\hat{X}_n$ and $X_n \hat{+} X_{n-2}$. Let $s_n$ be the length of the longest side for $\triangle(a_0, \ldots, a_n)$. If the triangle sequence uniquely describes a point, then $\lim_{n \to \infty} s_n = 0$. If the triangle sequence does not uniquely describe a point but instead describes a line segment $L$, of length, say, $l$, then we have $\lim_{n \to \infty} s_n = l$. We want to show that the even vertices $\hat{X}_{2n}$ converge to a point and that the odd vertices $\hat{X}_{2n+1}$ converge to a point, and further that each converges to one of the endpoints of the segment $L$. This will take some work.

**Lemma 4** *For all $n$, the point $X_n \hat{+} X_{n+2}$ is closer to the point $\hat{X}_{n+2}$ than to the point $\hat{X}_n$.*

The idea is that $X_n \hat{+} X_{n+2}$ is a weighted average of the vectors $\hat{X}_{n+2}$ and $\hat{X}_n$. Since $X_{n+2}$ is a longer vector than $X_n$, the result should be true. The actual proof is a straightforward calculation.

**Proof:** Denote the distance from a vector $X$ to a vector $Y$ by $d(X, Y)$.

By direct calculation, the vector from $\hat{X}_n$ to $X_n \hat{+} X_{n+2}$ is:

$$\begin{aligned}
\frac{1}{x_n + x_{n+2}}(X_{n+2} + X_n) - \frac{1}{x_n}X_n &= \frac{1}{x_n + x_{n+2}}X_{n+2} - \frac{x_{n+2}}{x_n(x_n + x_{n+2})}X_n \\
&= \frac{x_{n+2}}{x_n + x_{n+2}}\hat{X}_{n+2} - \frac{x_{n+2}}{x_n + x_{n+2}}\hat{X}_n \\
&= \frac{x_{n+2}}{x_n + x_{n+2}}(\hat{X}_{n+2} - \hat{X}_n),
\end{aligned}$$

and thus

$$d(\hat{X}_n, X_n \hat{+} X_{n+2}) = \frac{x_{n+2}}{x_n + x_{n+2}}d(\hat{X}_n, \hat{X}_{n+2}).$$

By a similar calculation, we have

$$d(\hat{X}_{n+2}, X_n \hat{+} X_{n+2}) = \frac{x_n}{x_n + x_{n+2}}d(\hat{X}_n, \hat{X}_{n+2}).$$

Since the $x_k$ are an increasing sequence, we have our result. $\square$

Our next lemma, whose proof we omit, is straightforward and is a simple geometric fact, but one which we will critically need.

**Lemma 5** *Let A,B and C be the three vertices of a triangle and let D be any point on the edge connecting the vertices B and C. Then*

$$d(A, D) \le \max(d(A, B), d(A, C)).$$

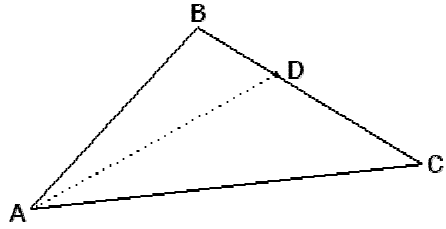

The theorem for this subsection is:

**Theorem 6**

$$\lim_{n \to \infty} d(\hat{X}_n, \hat{X}_{n+2}) = 0.$$

Note that this theorem is indeed simply stating that points $\hat{X}_k$ of the same parity converge.

**Proof:** Consider our triangle $\triangle(a_0, \dots, a_n)$.

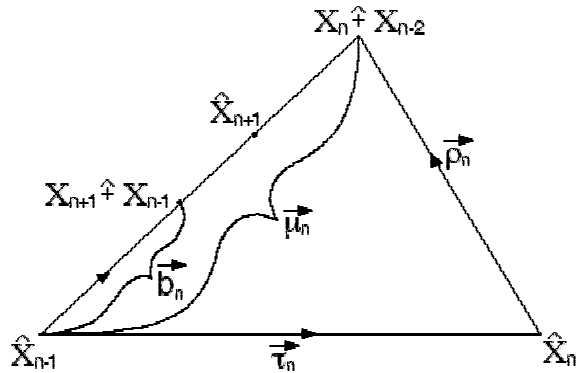

Set:

$$\vec{\rho}_n = \text{vector from } \hat{X}_n \text{ to } X_n \hat{+} X_{n-2}$$

$$\vec{\tau}_n = \text{vector from } \hat{X}_{n-1} \text{ to } \hat{X}_n$$

$$\vec{\mu}_n = \text{vector from } \hat{X}_{n-1} \text{ to } X_n \hat{+} X_{n-2}$$

$$\vec{b}_n = \text{vector from } \hat{X}_{n-1} \text{ to } X_{n-1} \hat{+} X_{n+1},$$

By $\rho_n$, we mean the length of the vector $\vec{\rho}_n$, etc.

We know from the first lemma of this subsection that $d(\hat{X}_n, X_n \hat{+} X_{n-2}) \leq d(\hat{X}_{n-2}, X_n \hat{+} X_{n-2})$. Then, since the points $\hat{X}_n, \hat{X}_{n-2}$ and $X_n \hat{+} X_{n-2}$ are collinear, we have that

$$\rho_n \leq \frac{1}{2} d(\hat{X}_n, \hat{X}_{n-2}).$$

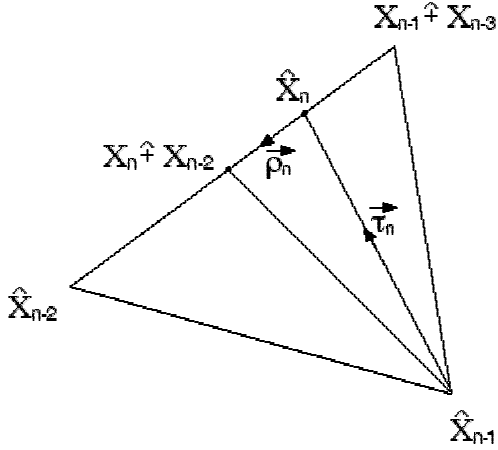

We have that the longest side lengths of each triangle, denoted by $s_n$, must approach $\ell$. If $\ell = 0$, then all of the triangles converge to a point and the lemma is true. Suppose, then, that $\ell \neq 0$. For any positve $\epsilon$, we can find an $N$ such that for all $n \geq N$,

$$\ell \leq s_n < \ell + \epsilon.$$

Choose any such $\epsilon$ such that $\epsilon < \ell$. From the above diagram we see that, for $n \geq N + 1$,

$$
\begin{aligned}
\rho_n &\leq \frac{1}{2}d(\hat{X}_n, \hat{X}_{n-2}) \\
&\leq \frac{1}{2}d(X_{n-1}\hat{+}X_{n-3}, \hat{X}_{n-2}) \\
&\leq \frac{1}{2}s_{n-1} \\
&\leq \frac{1}{2}(\ell + \epsilon) \\
&< \ell.
\end{aligned}
$$

Thus for large enough $n$, we have $\rho_n < \ell$ and $\rho_{n+1} < \ell$. This combined with the fact that $\ell \leq s_{n+1} = \max\{\tau_{n+1}, \mu_{n+1}, \rho_{n+1}\}$ shows that $\ell \leq \max\{\tau_{n+1}, \mu_{n+1}\}$. By our lemma on the triangle with vertices $A$,$B$ and $C$ (here the vertex $\hat{X}_n$ is playing the role of $A$ and the vertices $\hat{X}_{n-1}$ and $X_n \hat{+} X_{n-2}$ are playing the roles of $B$ and $C$), we have $\ell \leq \max\{\tau_n, \rho_n\}$. But $\rho_n < \ell$, meaning that $\tau_n \geq \ell$.

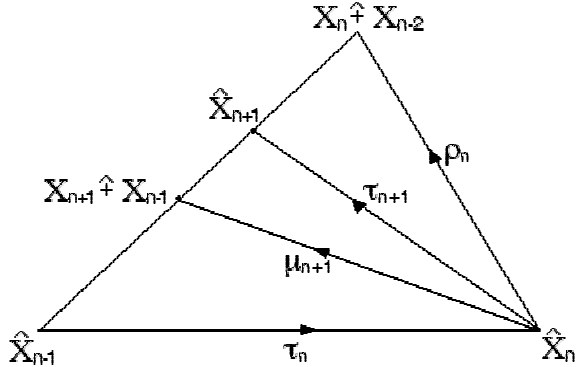

Assume for a moment that we can show, for large enough $n$, that $\tau_n \leq \mu_n$. Then for these large $n$, we know both that $s_n = \mu_n$ and that $\ell \leq \tau_n \leq s_n$. Then

$$\tau_n \to \ell.$$

Since the intersection of all of the $\triangle(a_0, \ldots, a_n)$ is the line segment $\ell$, we have our result, again provided that $\tau_n \leq \mu_n$. Thus we must prove this last inequality.

If the angle at the vertex $X_n \hat{+} X_{n-2}$ is obtuse or right, then we can see from a diagram similar to the one above that $\tau_{n+1} \leq \mu_{n+1}$.

Assume then that this angle is acute. Let $p$ be the foot of the perpendicular drawn from the point $\hat{X}_n$ to the line spanned by $\hat{X}_{n-1}$ and $X_n \hat{+} X_{n-2}$.

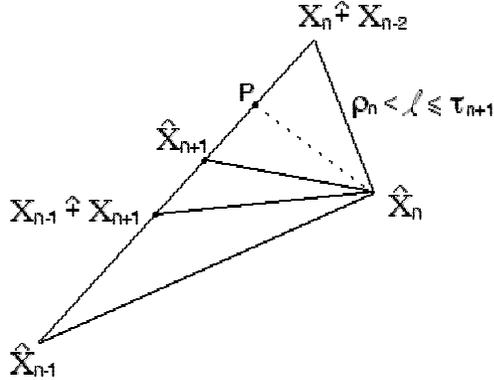

We have that the the point $\hat{X}_{n+1}$ is between $p$ and $X_{n-1} \hat{+} X_{n+1}$, since we know that $\rho_n < \tau_{n+1}$, giving us that $\tau_{n+1} \leq \mu_{n+1}$ is true also in this case. $\square$

## 6.2 First lemma towards uniqueness results

We will need the fact that making any finite number of changes in a triangle sequence will not effect questions of uniqueness. More precisely, we have the following:

**Lemma 7** *A triangle sequence $(a_0, a_1, a_2, \ldots)$ uniquely defines a pair of numbers if and only if the triangle sequence $(a_n, a_{n+1}, \ldots)$, for any $n > 0$, also uniquely describes a pair of numbers.*

This follows from the fact that locally, in the interior of any $\triangle_n$, the triangle map $T$ is bijective.

## 6.3 Uniqueness when $a_n = 0$ for infinitely many $n$

**Lemma 8** *Let $(a_0, a_1, a_2, \ldots)$ be a triangle sequence. If, for infinitely many of the $n$, we have $a_n = 0$, then the triangle sequence will describe a unique point.*

**Proof:** Recall that $s_n$ denotes the length of the longest edge of the triangle $\triangle(a_0, \ldots, a_n)$. We have seen that our triangle sequence will describe a unique point precisely when $\lim_{n \to \infty} s_n = 0$. Suppose that this does not happen. We keep the notation that $\ell = \lim_{n \to \infty} s_n$.

Set $\epsilon = \ell/2$. Recall that $\rho_n$ denotes the side length from the vertex $\hat{X}_n$ to $X_n \hat{+} X_{n-2}$ and that we have shown in the proof about the convergence of vertices of the same parity that

$$\lim_{n \to \infty} \rho_n = 0.$$

Then there exists a positive integer $M$, which we can make as large as we want, such that

$$
\begin{aligned}
\rho_M &< \epsilon \\
\rho_{M+1} &< \ell \\
a_{M+1} &= 0 \\
s_M &< \ell + \epsilon.
\end{aligned}
$$

Note that it is here that we are using our assumption that infinitely many of the $a_n$ are zero.

Since we always have that $X_{n+1} = X_{n-2} + a_{n+1}X_{n-1} + X_n$, we have

$$
X_{M+1} = X_{M-2} + X_M
$$

and thus $\tau_{M+1} = \rho_M$, where, recall, $\tau_{M+1}$ denotes the side length from the vertex $\hat{X}_M$ to the vertex $X_{M+1}\hat{+}X_{M-1}$. In the parity proof, we showed, for large enough $n$, that $\tau_n \geq \ell$. But then, choosing $M$ large enough, we have the desired contradiction

$$
\ell \leq \tau_{M+1} = \rho_M < \epsilon < \ell. \ \square
$$

## 6.4 $\Pi(1 - \lambda_n) = 0$ implies uniqueness

From the last two sections, we can assume that we have an infinite triangle sequence $(a_0, \dots)$ such that $a_n \neq 0$ for all $n$. Define for each $\triangle(a_0, \dots, a_n)$ the ratio

$$
\lambda_n = \frac{\text{Distance from } \hat{X}_{n-1} \text{ to } \hat{X}_{n+1}}{\text{Distance from } \hat{X}_{n-1} \text{ to } X_n\hat{+}X_{n-2}}.
$$

We will see that the question of uniqueness is linked to the size of the various $\lambda_n$. The goal of this section and the next is:

**Theorem 9** *Assume that $(a_0, \ldots)$ is a triangle sequence such that for all $n$, $a_n \neq 0$. Then this triangle sequence describes a unique pair $(\alpha, \beta)$ precisely when*

$$\prod_{n=0}^{\infty}(1 - \lambda_n) = 0.$$

In this section we show that if $\prod(1 - \lambda_n) = 0$, then we have uniqueness. In the next section we show that if $\prod(1 - \lambda_n) \neq 0$, then we have non-uniqueness. Then we will show that the infinite product $\prod(1 - \lambda_n) \neq 0$ when the individual $a_n$ grow sufficiently fast.

For the rest of this section, assume that $\prod(1 - \lambda_n) = 0$. We continue to use the notation that $\tau_n$ is the length of the side from the vertex $\hat{X}_{n-1}$ to the vertex $\hat{X}_n$. We have shown that

$$\lim_{n \to \infty} \tau_n = \ell$$

and that uniqueness is equivalent to $\ell = 0$.

We will break the proof into a number of lemmas involving inequalities. Assume for a moment the following lemma:

**Lemma 10** *For large enough $n$, assume that*

$$\frac{\tau_{n+1}}{\tau_n} \leq \sqrt{1 - \lambda_n}.$$

*Then*

$$\lim_{n \to \infty} \tau_n = 0,$$

*and thus the triangle sequence will describe a unique pair $(\alpha, \beta)$.*

Assuming this lemma, we have for any fixed $M$,

$$
\begin{aligned}
\tau_n &= \tau_M \prod_{i=M}^{n-1} \frac{\tau_{i+1}}{\tau_i} \\
&< \tau_M \sqrt{\prod_{i=M}^{n-1}(1-\lambda_i)} \\
&= \tau_M \left(\prod_{i=0}^{M-1}(1-\lambda_i)\right)^{-\frac{1}{2}} \sqrt{\prod_{i=0}^{n-1}(1-\lambda_i)},
\end{aligned}
$$

which converges to 0 as $n \to \infty$. (Here we used the fact that $a_i > 0$ for all $i$, so that $1 - \lambda_i > 0$ for all $i$.) Thus $\lim_{n\to\infty} \tau_n = 0$, which in turn means that $\ell = 0$ and that the triangle sequence $\{a_n\}$ corresponds to a unique point.

To prove that $\frac{\tau_{n+1}}{\tau_n} \leq \sqrt{1-\lambda_n}$, we need:

**Lemma 11** *If, for all $n$, we have*

$$
1 - \lambda_n \left(1 - \frac{\rho_n}{\tau_n}\right) \leq \sqrt{1-\lambda_n},
$$

*then*

$$
\frac{\tau_{n+1}}{\tau_n} \leq \sqrt{1-\lambda_n}.
$$

**Proof:** This will be a simple geometric argument using $\triangle(a_0, \ldots, a_n)$. Our notation is such that $\tau_{n+1}$ is the length of the vector $\vec{\tau}_{n+1}$ from the point $\hat{X}_n$ to the point $\hat{X}_{n+1}$.

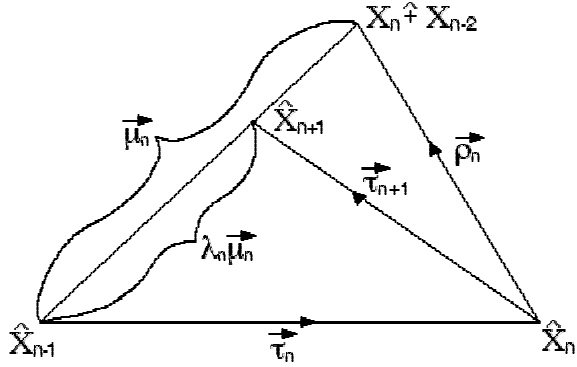

Then

$$\vec{\tau}_{n+1} = \lambda_n \vec{\mu}_n - \vec{\tau}_n.$$

But $\vec{\mu}_n = \vec{\tau}_n + \vec{\rho}_n$. Then we have

$$\vec{\tau}_{n+1} = \lambda_n \vec{\rho}_n - (1 - \lambda_n)\vec{\tau}_n.$$

Then

$$
\begin{aligned}
\frac{\tau_{n+1}}{\tau_n} &= \frac{1}{\tau_n}\left|\lambda_n \vec{\rho}_n - (1 - \lambda_n)\vec{\tau}_n\right| \\
&\leq \frac{1}{\tau_n}(\lambda_n \rho_n + (1 - \lambda_n)\tau_n) \\
&= 1 - \lambda_n\left(1 - \frac{\rho_n}{\tau_n}\right) \\
&< \sqrt{1 - \lambda_n}.
\end{aligned}
$$

and we are done with the lemma.

Thus we need to show that for large enough $n$, $1 - \lambda_n\left(1 - \frac{\rho_n}{\tau_n}\right) \leq \sqrt{1 - \lambda_n}$. This will take some work.

If we somehow know that $\ell = 0$, we already know that the triangle sequence describes a unique point. We can assume, then, that $\ell > 0$. Since $\lim_{n\to\infty} \rho_n = 0$, there is some $M$ such that $\rho_n < \frac{\ell}{3}$ for all $n \geq M$. Let $P$ be the point closest to $\hat{X}_{n-1}$ on the ray $\hat{X}_{n-1}(X_n\hat{+}X_{n-2})$ such that $d(\hat{X}_n, P) = 2\rho_n$.

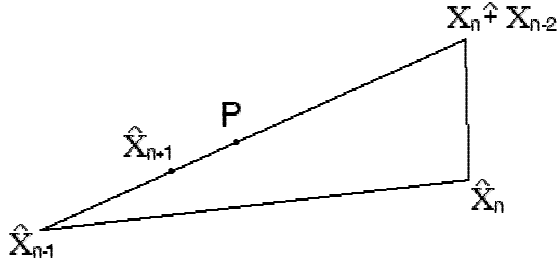

Since

$$d(\hat{X}_n, \hat{X}_{n-1}) = \tau_n > \ell > 2\rho_n = d(\hat{X}_n, P) > \rho_n = d(\hat{X}_n, X_n\hat{+}X_{n-2}),$$

we know that $P$ is on the line segment $\overline{\hat{X}_{n-1}(X_n\hat{+}X_{n-2})}$.

By choosing $M$ large enough, we can see that $\hat{X}_{n+1}$ is on the segment $\overline{\hat{X}_{n-1}P}$. Thus

$$
\begin{aligned}
(1-\lambda_n)\mu_n &= d(\hat{X}_{n+1}, X_n\hat{+}X_{n-2}) \\
&> d(P, X_n\hat{+}X_{n-2}) \\
&\geq d(P, \hat{X}_n) - d(X_n\hat{+}X_{n-2}, \hat{X}_n) \\
&= 2\rho_n - \rho_n = \rho_n.
\end{aligned}
$$

Since $\mu_n = d(\hat{X}_{n-1}, X_n\hat{+}X_{n-2}) \leq d(\hat{X}_{n-1}, \hat{X}_n) + d(X_n\hat{+}X_{n-2}, \hat{X}_n) = \tau_n + \rho_n$, we have $1 - \lambda_n > \frac{\rho_n}{\tau_n + \rho_n}$ and hence

$$\lambda_n < 1 - \frac{\rho_n}{\tau_n + \rho_n}.$$

We claim that $1 - \lambda_n \left(1 - \frac{\rho_n}{\tau_n}\right) < \sqrt{1 - \lambda_n}$, which is the inequality that we need to finish the proof of the theorem. The claim is equivalent to

$$1 - 2\lambda_n \left(1 - \frac{\rho_n}{\tau_n}\right) + \lambda_n{}^2 \left(1 - \frac{\rho_n}{\tau_n}\right)^2 < 1 - \lambda_n$$

$$\Longleftrightarrow \qquad \lambda_n{}^2 \left(1 - \frac{\rho_n}{\tau_n}\right)^2 < \lambda_n \left(1 - \frac{2\rho_n}{\tau_n}\right)$$

$$\Longleftrightarrow \qquad \lambda_n < \frac{1 - \frac{2\rho_n}{\tau_n}}{1 - \frac{2\rho_n}{\tau_n} + \frac{\rho_n{}^2}{\tau_n{}^2}} = 1 - \frac{\rho_n{}^2}{(\tau_n - \rho_n)^2}$$

Thus it suffices to show that $1 - \frac{\rho_n}{\tau_n + \rho_n} \leq 1 - \frac{\rho_n{}^2}{(\tau_n - \rho_n)^2}$. But that is equivalent to

$$\frac{\rho_n{}^2}{(\tau_n - \rho_n)^2} \leq \frac{\rho_n}{\tau_n + \rho_n}$$

$$\Longleftrightarrow \qquad \rho_n(\tau_n + \rho_n) \leq (\tau_n - \rho_n)^2$$

$$\Longleftrightarrow \quad 0 \leq \tau_n{}^2 - 3\tau_n \rho_n = \tau_n(\tau_n - 3\rho_n),$$

which is true because $\tau_n - 3\rho_n > \ell - 3 \cdot \frac{\ell}{3} = 0$. Our claim and hence theorem is proved.

## 6.5 $\Pi(1 - \lambda_n) \neq 0$ implies non-uniqueness

This is the most complicated section of this paper. Our goal is:

**Theorem 12** *Suppose that the sequence $\{a_n\}$ contains at most a finite number of zeros, such that $a_n > 0$ for $n > N$ and that*

$$\prod_{n=N}^{\infty} (1 - \lambda_n) > 0.$$

*Then the triangle sequence $\{a_n\}$ does not correspond to a unique point.*

As seen earlier, we can assume that $a_n \neq 0$ for all $n$. We will show nonuniqueness by showing, under the hypothesis of the theorem, that

$$\lim_{n \to \infty} \tau_n = \ell > 0.$$

As in the proof of the converse, this argument will come down to finding bounds on the ratios $\frac{\tau_{n+1}}{\tau_n}$, so that we will be able to reduce the above theorem to the following lemma:

**Lemma 13** *If for large $n$, we have the bounds*

$$\frac{\tau_{n+1}}{\tau_n} > (1 - \lambda_n)^2,$$

*then*

$$\lim_{n \to \infty} \tau_n = \ell > 0.$$

**Proof of lemma:** For large enough $n$, we have

$$
\begin{aligned}
\tau_n &= \tau_M \prod_{i=M}^{n-1} \frac{\tau_{i+1}}{\tau_i} \\
&> \tau_M \prod_{i=M}^{n-1} (1 - \lambda_i)^2 \\
&= \tau_M \left( \prod_{i=0}^{M-1} (1 - \lambda_i) \right)^{-2} \left( \prod_{i=0}^{n-1} (1 - \lambda_i) \right)^2,
\end{aligned}
$$

which converges to a positive constant. Hence $\ell = \lim_{n \to \infty} \tau_n$ is a positive constant, and we have non-uniqueness. $\square$

It will take serious work to show that $\frac{\tau_{n+1}}{\tau_n} > (1 - \lambda_n)^2$ for large $n$. First, if $\prod_{n=0}^{\infty} (1 - \lambda_n) > 0$, then $1 - \lambda_n \leq \frac{3}{4}$ for finitely many values of $n$. Thus there exists an $M$ such that $1 - \lambda_n > \frac{3}{4}$ whenever $n \geq M$. Thus $\frac{1}{\lambda_n} > 4$ whenever $n \geq M$.

We will show first:

**Lemma 14** *There exists an $M' \geq M$ such that*

$$\rho_{M'} \leq 2\tau_{M'}.$$

**Proof of lemma:** The proof is not at all obvious, but the heart of it lies in distinguishing six cases and dealing with each geometrically.

Set

$$\gamma_n = \angle \hat{X}_{n-1}\hat{X}_n(X_n \hat{+} X_{n-2})$$

$$\phi_n = \angle \hat{X}_n\hat{X}_{n-1}(X_n \hat{+} X_{n-2}),$$

$$\phi'_n = \angle \hat{X}_n(X_n \hat{+} X_{n-2})\hat{X}_{n-1}.$$

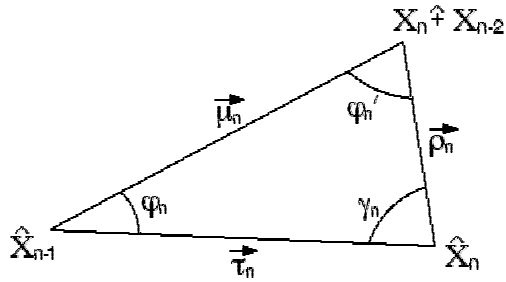

Let $P_n$ be the foot of the perpendicular from $\hat{X}_n$ to line $\overline{\hat{X}_{n-1}(X_n \hat{+} X_{n-2})}$.

The relationships between different lengths are configuration-dependent, so we will need to treat the different configurations as the six separate cases listed below. Intuitively, the worst cast scenario, Case 6, is when $\angle P_M \hat{X}_M \hat{X}_{M-1}$ is large, allowing $\tau_{M+1}$ to approximate the tiny height $P_M \hat{X}_M$ while leaving $\rho_{M+1}$ to be possibly the same order of magnitude as the long $P_M \hat{X}_{M+1}$. Fortunately, the other cases are not difficult, and Case 6 eventually stops occurring after a finite number of steps.

1. $\rho_M \leq 2\tau_M$.

   Then let $M' = M$.

2. $\rho_M > 2\tau_M$ and $\phi_M \geq \frac{\pi}{2}$.

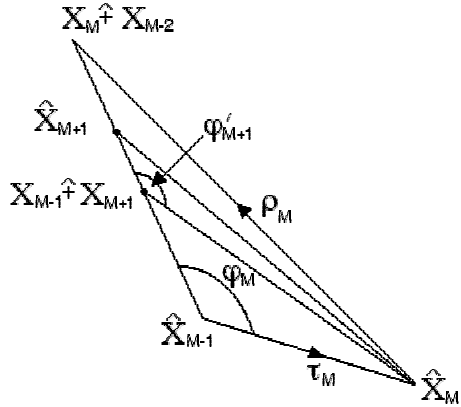

   The latter condition implies that $\phi'_{M+1}$ is obtuse, so that $\overline{\hat{X}_M \hat{X}_{M+1}}$ is the longest side of $\triangle \hat{X}_M \hat{X}_{M+1}(X_{M-1}\hat{+}X_{M+1})$. This means that $\rho_{M+1} < \tau_{M+1} < 2\tau_{M+1}$. Thus we have case 1 for $M+1$ and can let $M' = M+1$.

3. $\rho_M > 2\tau_M$, $\phi_M < \frac{\pi}{2}$, and $\gamma_M \leq \frac{3\pi}{4}$.

In the following, we will be freely using the inequalities shown in the appendix. Note that $\rho_M > 2\tau_M$ implies that $\tau_M$ is not the longest side of $\triangle \hat{X}_{M-1}\hat{X}_M(X_M\hat{+}X_{M-2})$. Therefore, $\phi'_M < \frac{\pi}{2}$ and $P_M$ lies on line segment $\overline{\hat{X}_{M-1}(X_M\hat{+}X_{M-2})}$. (The same will be true for the remaining cases.)

Let $y = d(\hat{X}_M, P_M)$. Then $\rho_M > \tau_M$ means that

$$\cos \angle P_M\hat{X}_M(X_M\hat{+}X_{M-2}) = \frac{y}{\rho_M} < \frac{y}{\tau_M} = \cos \angle P_M\hat{X}_M\hat{X}_{M-1}.$$

Since cosine is a decreasing function from 0 to $\pi/2$, we have

$$\angle P_M\hat{X}_M(X_M\hat{+}X_{M-2}) > \angle P_M\hat{X}_M\hat{X}_{M-1}.$$

Then

$$\angle P_M\hat{X}_M\hat{X}_{M-1} < \frac{1}{2}(\angle P_M\hat{X}_M(X_M\hat{+}X_{M-2}) + \angle P_M\hat{X}_M\hat{X}_{M-1}),$$

which in turn, since $(\angle P_M\hat{X}_M(X_M\hat{+}X_{M-2}) + \angle P_M\hat{X}_M\hat{X}_{M-1}) = \gamma_M$,

$$\angle P_M\hat{X}_M\hat{X}_{M-1} < \frac{1}{2}\gamma_M \leq \frac{3\pi}{8}.$$

This gives us the bound that we will need in the next paragraph:

$$\tau_M = \frac{y}{\cos \angle P_M\hat{X}_M\hat{X}_{M-1}} \leq \frac{y}{\cos \frac{3\pi}{8}} \leq \frac{\tau_{M+1}}{\cos \frac{3\pi}{8}}.$$

Using lemma 4 from section 6.1, we have:

$$\begin{aligned}
\rho_{M+1} &\leq \frac{1}{2}d(\hat{X}_{M+1}, \hat{X}_{M-1}) \\
&\leq \frac{1}{2}(\tau_{M+1} + \tau_M) \\
&\leq \frac{1}{2}\left(\tau_{M+1} + \frac{\tau_{M+1}}{\cos \frac{3\pi}{8}}\right) \\
&< 2\tau_{M+1}.
\end{aligned}$$

Thus we have case 1 for $M + 1$ and can let $M' = M + 1$.

4. $\rho_M > 2\tau_M$, $\phi_M < \frac{\pi}{2}$, $\gamma_M > \frac{3\pi}{4}$, and $X_{M-1}\hat{+}X_{M+1}$ lies on line segment $\overline{(X_M\hat{+}X_{M-2})P_M}$.

   Then $\phi'_{M+1} = \angle\hat{X}_M(X_{M-1}\hat{+}X_{M+1})\hat{X}_{M+1}$ is obtuse, making $\overline{\hat{X}_M\hat{X}_{M+1}}$ the longest side of $\triangle\hat{X}_M\hat{X}_{M+1}(X_{M-1}\hat{+}X_{M+1})$. This means that $\rho_{M+1} < \tau_{M+1} < 2\tau_{M+1}$. Thus we have case 1 for $M + 1$ and can simply let $M' = M + 1$.

5. $\rho_M > 2\tau_M$, $\phi_M < \frac{\pi}{2}$, $\gamma_M > \frac{3\pi}{4}$, $X_{M-1}\hat{+}X_{M+1}$ lies on line segment $\overline{P_M\hat{X}_{M-1}}$, and $\pi - \gamma_{M+1} \leq 2(\pi - \gamma_M)$.

   We will be freely using the numbers $\lambda_n$, $\tilde{\lambda}_n$ and $\lambda'_n$, which are defined in the appendix. Since $\rho_M > \tau_M$ and since $\phi_M < \frac{\pi}{2}$, we see that $\frac{\pi}{2} > \phi_M > \phi'_M$. Thus

   $$
   \begin{aligned}
   \frac{d(P_M, \hat{X}_{M-1})}{d(X_M\hat{+}X_{M-2}, P_M)} &= \frac{\tau_M}{\rho_M} \cdot \frac{d(P_M, \hat{X}_{M-1})/\tau_M}{d(X_M\hat{+}X_{M-2}, P_M)/\rho_M} \\
   &= \frac{\tau_M}{\rho_M} \cdot \frac{\cos\phi_M}{\cos\phi'_M} \\
   &< \frac{\tau_M}{\rho_M}.
   \end{aligned}
   $$

This implies that

$$
\begin{aligned}
\tilde{\lambda}_M &= \frac{d(X_{M-1}\hat{+}X_{M+1}, \hat{X}_{M-1})}{d(X_M\hat{+}X_{M-2}, \hat{X}_{M-1})} \\
&\leq \frac{d(P_M, \hat{X}_{M-1})}{d(X_M\hat{+}X_{M-2}, P_M) + d(P_M, \hat{X}_{M-1})} \\
&< \frac{d(P_M, \hat{X}_{M-1})}{\frac{\rho_M}{\tau_M}d(P_M, \hat{X}_{M-1}) + d(P_M, \hat{X}_{M-1})}
\end{aligned}
$$

$$= \frac{\tau_M}{\rho_M + \tau_M}.$$

From the definition of $\lambda'_M$ and from lemma 25, in the appendix, we have:

$$\rho_{M+1} = \lambda'_M \mu_M \leq \frac{\mu_M}{\frac{1}{\lambda_M}\left(\frac{1}{\lambda_M} - 1\right)}.$$

Since $\mu_M < \rho_M + \tau_M$ (these are the three side lengths of a triangle) and using the above inequality on $\tilde{\lambda}_M$, we have

$$\rho_{M+1} < \frac{\rho_M + \tau_M}{\frac{\rho_M + \tau_M}{\tau_M}\left(\frac{\rho_M + \tau_M}{\tau_M} - 1\right)} = \frac{\tau_M{}^2}{\rho_M}.$$

Using that $\frac{\sin\phi_M}{\tau_{M+1}} = \frac{\sin\gamma_{M+1}}{\tau_M}$ and that $\frac{\sin\gamma_M}{\mu_M} = \frac{\sin\phi_M}{\rho_M}$ (both following from the law of sines), we also have

$$\begin{aligned}
\tau_{M+1} &= \tau_{M+1}\left(\frac{\sin\phi_M}{\tau_{M+1}} \cdot \frac{\tau_M}{\sin\gamma_{M+1}}\right)\left(\frac{\rho_M}{\sin\phi_M} \cdot \frac{\sin\gamma_M}{\mu_M}\right) \\
&= \frac{\tau_M \rho_M}{\mu_M} \cdot \frac{\sin(\pi - \gamma_M)}{\sin(\pi - \gamma_{M+1})}.
\end{aligned}$$

But $0 < \pi - \gamma_{M+1} \leq 2(\pi - \gamma_M) < \frac{\pi}{2}$, so

$$\frac{\sin(\pi - \gamma_M)}{\sin(\pi - \gamma_{M+1})} \geq \frac{\sin(\pi - \gamma_M)}{\sin(2(\pi - \gamma_M))} = \frac{1}{2\cos(\pi - \gamma_M)} \geq \frac{1}{2}.$$

Hence

$$\tau_{M+1} \geq \frac{\tau_M \rho_M}{2\mu_M} > \frac{\tau_M \rho_M}{2(\rho_M + \tau_M)}.$$

Therefore,

$$\begin{aligned}
\frac{\rho_{M+1}}{\tau_{M+1}} &< \frac{\tau_M{}^2}{\rho_M} \cdot \frac{2(\rho_M + \tau_M)}{\tau_M \rho_M} = \frac{2\tau_M(\rho_M + \tau_M)}{\rho_M{}^2} \\
&< \frac{\rho_M\left(\rho_M + \frac{1}{2}\rho_M\right)}{\rho_M{}^2} = \frac{3}{2},
\end{aligned}$$

and $\rho_{M+1} < 2\tau_{M+1}$. Thus we have case 1 for $M+1$ and can let $M' = M+1$.

6. $\rho_M > 2\tau_M$, $\phi_M < \frac{\pi}{2}$, $\gamma_M > \frac{3\pi}{4}$, $X_{M-1}\hat{+}X_{M+1}$ lies on line segment $\overline{P_M\hat{X}_{M-1}}$, and $\pi - \gamma_{M+1} > 2(\pi - \gamma_M)$.

   We will see that this case cannot occur for all $n \geq M$. Suppose that it does. Then we will get $\pi - \gamma_{n+1} > 2(\pi - \gamma_n)$ for $n \geq M$. Thus, we have $\pi - \gamma_{M+i} > 2^i(\pi - \gamma_M)$. But then we can make $\pi - \gamma_{M+i}$ as large as we like, by choosing large enough $i$. This is not possible, implying that we cannot be in case six for all $n \geq M$. But then there is some $n \geq M$ such that we are in one of the first five cases, in which case we know we are done.

$\square$

Our next technical lemma is:

**Lemma 15** *Assume that $\rho_M \leq 2\tau_M$ and that, for all $n \geq M$, $\lambda_n \leq \frac{1}{4}$. Then for all $n \geq M + 2$, we have*

$$\rho_n \leq \frac{1}{2}\tau_n.$$

**Proof of lemma:**

We argue by induction on $n$. We will be freely using the inequalities of the last lemma in the second section of the appendix. For the base case, we first find the inequality for $\rho_{M+1} \leq \frac{21}{40}\tau_{M+1}$ and then show $\rho_{M+2} \leq \frac{1}{2}\tau_{M+2}$.. We have from lemma 26 in the appendix that

$$\tau_{M+1} \geq \tau_M - \lambda_M(\rho_M + \tau_M) \geq \tau_M - \frac{1}{4}(2\tau_M + \tau_M) = \frac{1}{4}\tau_M$$

and, also from lemma 26,

$$\rho_{M+1} \leq \frac{\rho_M + \tau_M}{\frac{1}{\lambda_M}\left(\frac{1}{\lambda_M}+1\right)} \leq \frac{2\tau_M + \tau_M}{4(4+1)} = \frac{3}{20}\tau_M.$$

Hence

$$\rho_{M+1} - \frac{1}{2}\tau_{M+1} \leq \left(\frac{3}{20} - \frac{1}{2}\cdot\frac{1}{4}\right)\tau_M = \frac{1}{40}\tau_M$$

Using that $\tau_M \leq 4\tau_{M+1}$, then

$$
\begin{aligned}
\rho_{M+1} &\leq \frac{1}{2}\tau_{M+1} + \frac{1}{40}\tau_M \\
&\leq (\frac{1}{2} + \frac{1}{10})\tau_{M+1} \\
&= \frac{3}{5}\tau_{M+1}.
\end{aligned}
$$

We have

$$\tau_{M+2} \geq \tau_{M+1} - \lambda_{M+1}(\rho_{M+1} + \tau_{M+1}) \geq \tau_{M+1} - \frac{1}{4}\left(\frac{3}{5}\tau_{M+1} + \tau_{M+1}\right) = \frac{3}{5}\tau_{M+1}$$

and

$$\rho_{M+2} \leq \frac{\rho_{M+1} + \tau_{M+1}}{\frac{1}{\lambda_{M+1}}\left(\frac{1}{\lambda_{M+1}}+1\right)} \leq \frac{\frac{3}{5}\tau_{M+1} + \tau_{M+1}}{4(4+1)} = \frac{8}{100}\tau_{M+1}.$$

Then

$$\rho_{M+2} - \frac{1}{2}\tau_{M+2} \leq \left(\frac{8}{100} - \frac{1}{2}\cdot\frac{3}{5}\right)\tau_{M+1} < 0$$

and the base case is proved.

Now suppose that $\rho_i \leq \frac{1}{2}\tau_i$ for some $i \geq M+2$. Then

$$\tau_{i+1} \geq \tau_i - \lambda_i(\rho_i + \tau_i) \geq \tau_i - \frac{1}{4}\left(\frac{1}{2}\tau_i + \tau_i\right) = \frac{5}{8}\tau_i$$

and

$$\rho_{i+1} \leq \frac{\rho_i + \tau_i}{\frac{1}{\lambda_i}\left(\frac{1}{\lambda_i}+1\right)} \leq \frac{\frac{1}{2}\tau_i + \tau_i}{4(4+1)} = \frac{3}{40}\tau_i.$$

Hence

$$\rho_{i+1} - \frac{1}{2}\tau_{i+1} \leq \left(\frac{3}{40} - \frac{1}{2} \cdot \frac{5}{8}\right)\tau_i = -\frac{19}{80}\tau_i < 0$$

and induction is complete. $\square$

**Lemma 16** *With the same assumptions as in the previous lemma, we have*

$$1 - \lambda_n\left(1 + \frac{\rho_n}{\tau_n}\right) > (1 - \lambda_n)^2$$

*for $n \geq M' + 1$.*

**Proof of lemma:**

This equality is equivalent to

$$1 - \lambda_n\left(1 + \frac{\rho_n}{\tau_n}\right) > 1 - 2\lambda_n + \lambda_n{}^2$$

$$\Longleftrightarrow \qquad \lambda_n{}^2 < \lambda_n\left(1 - \frac{\rho_n}{\tau_n}\right)$$

$$\Longleftrightarrow \qquad \lambda_n < 1 - \frac{\rho_n}{\tau_n},$$

which is true since

$$\lambda_n < \frac{1}{4} < 1 - \frac{1}{2} \leq 1 - \frac{\rho_n}{\tau_n}.$$

$\square$

The above lemma is important since for $n \geq M' + 1$, by Lemma 26 in the appendix,

$$\frac{\tau_{n+1}}{\tau_n} \geq 1 - \lambda_n\left(1 + \frac{\rho_n}{\tau_n}\right) > (1 - \lambda_n)^2.$$

which is precisely the inequality needed to show non-uniqueness.

## 6.6 Explicit Examples of Non-uniqueness

The following restatement sums up the above results:

The triangle sequence $\{a_n\}$ does not correspond to a unique point if and only if it contains a finite number of zeros and $\prod_{n=N}^{\infty}(1 - \lambda_n) > 0$ (where $N$ is such that $a_n > 0$ for $n > N$).

But what is this nebulous $1 - \lambda_n$ thing? There turns out to be a nice simplification of this expression:

$$
\begin{aligned}
1 - \lambda_n &= 1 - \frac{\frac{x_n + x_{n-2}}{x_{n-1}}}{a_{n+1} + \frac{x_n + x_{n-2}}{x_{n-1}}} \\
&= \frac{a_{n+1}}{a_{n+1} + \frac{x_n + x_{n-2}}{x_{n-1}}} \\
&= \frac{a_{n+1} x_{n-1}}{a_{n+1} x_{n-1} + x_n + x_{n-2}} \\
&= \frac{a_{n+1} x_{n-1}}{x_{n+1}},
\end{aligned}
$$

since, by Proposition 1, we always have $x_{n+1} = a_{n+1} x_{n-1} + x_n + x_{n-2}$. Hence

$$
\begin{aligned}
\prod_{n=N}^{\infty} (1 - \lambda_n) &= \lim_{M \to \infty} \prod_{n=N}^{M} \frac{a_{n+1} x_{n-1}}{x_{n+1}} \\
&= \lim_{M \to \infty} \frac{x_{N-1} x_N}{x_M x_{M+1}} \prod_{n=N+1}^{M+1} a_n \\
&= x_{N-1} x_N \lim_{M \to \infty} \frac{1}{x_M x_{M+1}} \prod_{n=N+1}^{M+1} a_n.
\end{aligned}
$$

The question of uniqueness thus boils down to whether or not this limit is zero. For the following examples we obtain estimates on $1 - \lambda_n$ in order to use the above criterion. Again, these are the first examples of non-uniqueness found, though other non-unique triangle sequences are easy to generate empirically using the Mathematica package at:

http://www.williams.edu/Mathematics/tgarrity/triangle.html

**6.6.1** $a_n = n$

We claim that the triangle sequence $\{1, 2, 3, \ldots\}$ corresponds to a unique point. Observe that

$$1 - \lambda_n = \frac{a_{n+1}}{a_{n+1} + \frac{x_n + x_{n-2}}{x_{n-1}}} \leq \frac{a_{n+1}}{a_{n+1} + 1} = \frac{n+1}{n+2}.$$

(Here we are using that the $x_n$ form an increasing sequence of integers.) Thus

$$\prod_{n=0}^{\infty} (1 - \lambda_n) \leq \frac{1}{2} \cdot \frac{2}{3} \cdot \frac{3}{4} \cdot \ldots = 0$$

**6.6.2** $a_n = n^2$

We claim that the triangle sequence $\{1, 4, 9, \ldots\}$ corresponds to a unique point.

We first must find a lower bound on various $\frac{x_n}{x_{n-1}}$. More precisely, as we will see, we need to show that, for each $n$, either

$$\frac{x_n}{x_{n-1}} \geq n + 1$$

or

$$\frac{x_{n-1}}{x_{n-2}} \geq n.$$

This follows from

$$\begin{aligned}
\frac{x_{n-1}}{x_{n-2}} + \frac{x_n}{x_{n-1}} &= \frac{x_{n-1}}{x_{n-2}} + \frac{x_{n-1} + a_n x_{n-2} + x_{n-3}}{x_{n-1}} \\
&\geq \frac{x_{n-1}}{x_{n-2}} + 1 + \frac{a_n}{\frac{x_{n-1}}{x_{n-2}}},
\end{aligned}$$

using that $x_{n-3} > 0$. But this last term is greater than or equal to

$$2(\frac{x_{n-1}}{x_{n-2}} \cdot \frac{a_n}{\frac{x_{n-1}}{x_{n-2}}})^{1/2} + 1 = n + (n+1),$$

giving us our bound on either $\frac{x_n}{x_{n-1}}$ or $\frac{x_{n-1}}{x_{n-2}}$.

We need one more bound before we show uniqueness. We know that there are infinitely many $n$ such that $\frac{x_n}{x_{n-1}} \geq n+1$. Choose such an $n$. Then

$$
\begin{aligned}
1 - \lambda_n \quad &= \quad \frac{a_{n+1}}{a_{n+1} + \frac{x_n + x_{n-2}}{x_{n-1}}} \\
&\leq \quad \frac{a_{n+1}}{a_{n+1} + \frac{x_n}{x_{n-1}}} \\
&\leq \quad \frac{(n+1)^2}{(n+1)^2 + (n+1)} \\
&= \quad 1 - \frac{1}{n+2}.
\end{aligned}
$$

Now we know for each $k$, that either $\frac{x_{2k}}{x_{2k-1}} \geq 2k+1$ or $\frac{x_{2k-1}}{x_{2k-2}} \geq 2k$. Then, letting $n = 2k$ or $n = 2k-1$, from the the above inequality, we have

$$1 - \lambda_n \leq 1 - \frac{1}{2k+2}.$$

We want to show that

$$\prod_{n=0}^{\infty}(1 - \lambda_n) = 0.$$

No matter what is the integer $n$, we know that

$$1 - \lambda_n \leq 1.$$

Thus for $\prod_{n=0}^{\infty}(1 - \lambda_n)$, we will only take the product over those integers $n$ such that $\frac{x_n}{x_{n-1}} \geq n+1$. The terms that we have dropped can only make the product smaller.

Thus

$$
\begin{aligned}
\prod_{n=0}^{\infty} (1 - \lambda_n) \;&\leq\; \prod_{k=1}^{\infty} \left(1 - \frac{1}{2k+2}\right) \\
&=\; \frac{3}{4} \cdot \frac{5}{6} \cdot \frac{7}{8} \cdot \ldots \\
&\leq\; \left(\frac{3}{4} \cdot \frac{4}{5}\right)^{\frac{1}{2}} \cdot \left(\frac{5}{6} \cdot \frac{6}{7}\right)^{\frac{1}{2}} \cdot \left(\frac{7}{8} \cdot \frac{8}{9}\right)^{\frac{1}{2}} \cdot \ldots \\
&=\; 0.
\end{aligned}
$$

Thus this sequence corresponds to a unique point in the triangle.

### 6.6.3 $a_n = n$th prime

We claim that the triangle sequence consisting of the primes corresponds to a unique point. Let $p_n$ denote the $n$th prime. Observe that

$$
1 - \lambda_n = 1 - \frac{\frac{x_n + x_{n-2}}{x_{n-1}}}{a_{n+1} + \frac{x_n + x_{n-2}}{x_{n-1}}} \leq 1 - \frac{1}{a_{n+1} + 1} = 1 - \frac{1}{p_{n+1} + 1}.
$$

Further note that

$$
\begin{aligned}
\left(1 - \frac{1}{p}\right) - \left(1 - \frac{1}{p+1}\right)^2 \;&=\; \frac{(p+1)^2(p-1) - p^3}{p(p+1)^2} \\
&=\; \frac{p^2 - p - 1}{p(p+1)^2} \\
&>\; 0.
\end{aligned}
$$

Hence

$$
\prod_{n=0}^{\infty} (1 - \lambda_n) \leq \prod_{n=1}^{\infty} \left(1 - \frac{1}{p_n + 1}\right) \leq \left(\prod_{n=1}^{\infty} \left(1 - \frac{1}{p_n}\right)\right)^{\frac{1}{2}} = 0.
$$

Thus

$$
\prod_{n=0}^{\infty} (1 - \lambda_n) = 0.
$$

### 6.6.4   $a_n = 2^{n-1}$

Here we set $a_0 = 0$ and $a_n = 2^{n-1}$ for $n > 0$.

We claim that the triangle sequence $\{0, 1, 2, 4, 8, \ldots\}$ does not correspond to a unique point. We first establish by induction that for $n \geq 7$,

$$1 + 2^{\frac{n-1}{2}} \leq \frac{x_n}{x_{n-1}} \leq 2^{\frac{n}{2}}.$$

The base case $n = 7$ can be checked computationally as follows. By direct calculation we have

$$(x_0, x_1, x_2 \ldots) = (1, 1, 3, 8, 33, 164, 1228, 11757, \ldots),$$

and thus

$$1 + 2^{\frac{7-1}{2}} = 9 < \frac{11757}{1228} = \frac{x_7}{x_6} < 2^{\frac{7}{2}}.$$

Suppose our claim is true for some $k \geq 7$. Then

$$\frac{x_{k+1}}{x_k} = \frac{x_k + a_{k+1}x_{k-1} + x_{k-2}}{x_k} \geq 1 + \frac{a_{k+1}}{\frac{x_k}{x_{k-1}}} \geq 1 + \frac{2^k}{2^{\frac{k}{2}}} = 2^{\frac{k}{2}} + 1.$$

Also,

$$\frac{x_{k+1}}{x_k} \leq \frac{x_k + a_{k+1}x_{k-1} + x_{k-1}}{x_k} = 1 + \frac{a_{k+1} + 1}{\frac{x_k}{x_{k-1}}} \leq 1 + \frac{2^k + 1}{2^{\frac{k-1}{2}} + 1} \leq 2^{\frac{k+1}{2}}.$$

We have proven our estimates for $\frac{x_n}{x_{n-1}}$.

Using the upper bound for $\frac{x_n}{x_{n-1}}$, we can obtain a lower bound for $1 - \lambda_n$ for $n \geq 7$:

$$1 - \lambda_n \quad = \quad \frac{a_{n+1}}{a_{n+1} + \frac{x_n + x_{n-2}}{x_{n-1}}}$$

$$\geq \frac{a_{n+1}}{a_{n+1} + \frac{x_n + x_n}{x_{n-1}}}$$

$$\geq \frac{2^n}{2^n + 2 \cdot 2^{\frac{n}{2}}}$$

$$= \frac{2^n}{2^n + 2^{\frac{n}{2}+1}}$$

$$= 1 - \frac{1}{2^{\frac{n}{2}-1} + 1}$$

$$> 1 - 2^{-\frac{n}{2}+1}.$$

Hence

$$\prod_{n=7}^{\infty} (1 - \lambda_n) \geq \prod_{n=7}^{\infty} \left(1 - 2^{-\frac{n}{2}+1}\right) = \prod_{n=5}^{\infty} \left(1 - 2^{-\frac{n}{2}}\right).$$

Now

$$-\log\left(\prod_{n=5}^{\infty} \left(1 - 2^{-\frac{n}{2}}\right)\right) = \sum_{n=5}^{\infty} -\log\left(1 - 2^{-\frac{n}{2}}\right)$$

$$= \sum_{n \geq 5} \sum_{j \geq 1} \frac{2^{-\frac{jn}{2}}}{j}$$

$$= \sum_{j \geq 1} \frac{1}{j} \sum_{n \geq 5} 2^{-\frac{jn}{2}}$$

$$= \sum_{j \geq 1} \frac{1}{j} \cdot \frac{2^{-2j}}{1 - 2^{-\frac{j}{2}}}$$

$$\leq \sum_{j \geq 1} \frac{4}{j} \cdot \left(\frac{1}{4}\right)^j$$

$$= c$$

for a positive constant $c$. Therefore

$$\prod_{n=4}^{\infty} (1 - \lambda_n) \geq \prod_{n=4}^{\infty} \left(1 - 2^{-\frac{n}{2}}\right) \geq e^{-c} > 0.$$

# 7 Topological Dynamics of the Triangle Iteration

Since the traditional continued fraction algorithm gives an ergodic transformation of the unit interval, it is natural to ask about the dynamical properties of the triangle map. For most of the other types of multidimensional continued fractions, such questions have been asked and in fact answered. Although most of these algorithms have been shown to be ergodic (see F. Schweiger's work in [29], [28] and [26]), the techniques that are used do not appear to be immediately applicable to the triangle sequence, precisely because the analogue of uniqueness holds for these other algorithms. We are thus not yet able to determine whether or not the triangle sequence is ergodic, but can show that it is topologically strongly mixing, which implies that it is topologically ergodic and transitive. (We will give these definitions in a moment; a general reference is in [6], in chapter two, section 4.)

## 7.1 On open sets and partition triangles

This section will give us the needed lemmas for triangle sequences that will allow us to prove dynamical properties in the next section. We have partitioned our initial triangle $\triangle$ into infinitely many smaller triangles $\triangle(a_0, a_1, \ldots, a_n)$. We call these $\triangle(a_0, a_1, \ldots, a_n)$ *partition triangles.*

**Lemma 17** *The union of the edges of all partition triangles is dense in $\triangle$.*

**Proof:** The set of rational points is dense in $\triangle$. These points yield terminating triangle sequences [9] and are thus on the edges of partition triangles. $\square$

**Lemma 18** *The union of all partition triangle edges of the form*

$$\overline{\hat{X}_{n-1}(X_n\hat{+}X_{n-2})}$$

*is dense.*

**Proof:** Given an open ball $B_\epsilon(w)$ of radius $\epsilon$ about a point $w$ in the triangle, by the above lemma we know that there exists $w' \in B_\epsilon(w)$ that is on an edge of a partition triangle $\triangle\hat{X}_{n-1}\hat{X}_n(X_n\hat{+}X_{n-2})$ that corresponds to a terminating triangle sequence $\{a_0, \ldots, a_n\}$. From this triangle sequence, we want to construct a possibly new triangle sequence such that the appropiate edge of the new triangle intersects the open ball $B_\epsilon(w)$. The point $w'$ is on one of the three edges of the triangle formed from $\{a_0, \ldots, a_n\}$, giving us the following three cases to consider.

1. $w' \in \overline{\hat{X}_{n-1}(X_n\hat{+}X_{n-2})}$. Then we are done.

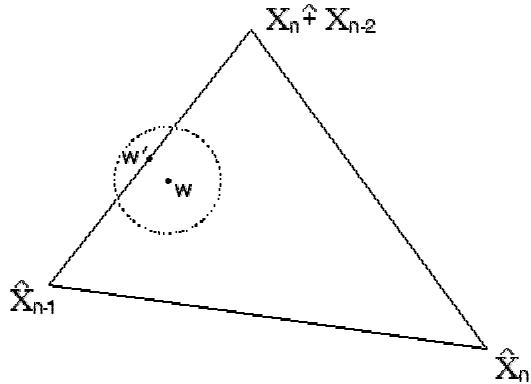

2. $w' \in \overline{\hat{X}_{n-1}\hat{X}_n}$. As $a_{n+1} \to \infty$, $X_{n-1}\hat{+}X_{n+1} \to \hat{X}_{n-1}$. But $\overline{\hat{X}_{n-1}\hat{X}_n}$ passes through $B_\epsilon(w)$ (i.e., at $w'$). Thus we can take a sufficiently large $a_{n+1}$ such that $\overline{(X_{n-1}\hat{+}X_{n+1})\hat{X}_n}$ also passes through $B_\epsilon(w)$.

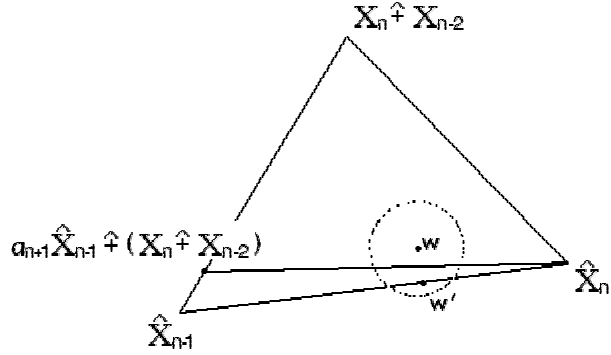

3. $w' \in \overline{\hat{X}_n(X_n\hat{+}X_{n-2})}$. Let $a_{n+1} = 0$. Then $X_n\hat{+}X_{n-2} = \hat{X}_{n+1}$, so that $w' \in \overline{\hat{X}_n\hat{X}_{n+1}}$. Hence we are back in case 2.

$\square$

**Theorem 19** *The set of all partition triangle vertices is dense in* $\triangle$*. In fact, the set of vertices of the form* $\hat{X}_n$ *are dense.*

**Proof:** Given any open ball $B_\epsilon(w)$, by the above lemma we know that there exists a triangle sequence $\{a_0, \ldots, a_n\}$ such that $\overline{\hat{X}_{n-1}(X_n\hat{+}X_{n-2})}$ intersects $B_{\frac{\epsilon}{2}}(w)$. Let $w_0 \in B_{\frac{\epsilon}{2}}(w) \cap \overline{\hat{X}_{n-1}(X_n\hat{+}X_{n-2})}$ and let $\ell_0 = d(\hat{X}_{n-1}, X_n\hat{+}X_{n-2})$.

We will inductively define the rest of the $a_i$'s, three at a time. Suppose that for some $k \geq 0$, we have defined $a_1, \ldots, a_{3k+n}$ and

$$w_k \in B_{\frac{\epsilon}{2}}(w) \cap \overline{\hat{X}_{3k+n-1}(X_{3k+n}\hat{+}X_{3k+n-2})}.$$

Then choose $a_{3k+n+1}$ such that

$$w_k \in \overline{\hat{X}_{3k+n+1}(X_{3k+n+1}\hat{+}X_{3k+n-1})}.$$

Next, let $a_{3k+n+2} = 0$, so that $X_{3k+n+1}\hat{+}X_{3k+n-1} = X_{3k+n+2}$ and $w_k \in \overline{\hat{X}_{3k+n+1}\hat{X}_{3k+n+2}}$. Lastly, choose a sufficiently large $a_{3k+n+3}$ so that

1. $B_{\frac{\epsilon}{2}}(w)$ and $\overline{\hat{X}_{3k+n+2}(X_{3k+n+3}\hat{+}X_{3k+n+1})}$ intersect (say at $w_{k+1}$), and

2. $d(\hat{X}_{3k+n+2}, X_{3k+n+3}\hat{+}X_{3k+n+1}) < \frac{3}{2}d(\hat{X}_{3k+n+2}, \hat{X}_{3k+n+1})$.

Our inductive definition is complete.

Let $\ell_k = d(\hat{X}_{3k+n-1}, X_{3k+n}\hat{+}X_{3k+n-2})$. Then

$$
\begin{aligned}
\ell_{k+1} &= d(\hat{X}_{3k+n+2}, X_{3k+n+3}\hat{+}X_{3k+n+1}) \\
&< \frac{3}{2}d(\hat{X}_{3k+n+2}, \hat{X}_{3k+n+1}) \\
&= \frac{3}{2}d(X_{3k+n+1}\hat{+}X_{3k+n-1}, \hat{X}_{3k+n+1}) \\
&\leq \frac{3}{2} \cdot \frac{1}{2}d(\hat{X}_{3k+n-1}, X_{3k+n}\hat{+}X_{3k+n-2}) \\
&= \frac{3}{4}\ell_k,
\end{aligned}
$$

using lemma 4 for the inequality $d(X_{3k+n+1}\hat{+}X_{3k+n-1}, \frac{1}{2}d(\hat{X}_{3k+n-1}, X_{3k+n}\hat{+}X_{3k+n-2})$. Hence $\ell_k < \left(\frac{3}{4}\right)^k \ell_0$, and thus $\ell_k \to 0$. Choose a large enough $k$ such that $\ell_k < \frac{\epsilon}{2}$. Then

$$
\begin{aligned}
d(X_{3k+n-1}, w) &\leq d(X_{3k+n-1}, w_k) + d(w_k, w) \\
&\leq \ell_k + \frac{\epsilon}{2} \\
&< \epsilon.
\end{aligned}
$$

Hence $B_\epsilon(w)$ contains the partition triangle vertex $X_{3k+n-1}$. $\square$

The key theorem for dynamical properties is:

**Theorem 20** *Any given open ball $B_\epsilon(w)$ contains a partition triangle.*

**Proof:** By the above theorem, we can choose $a_1, \ldots, a_n$ such that $\hat{X}_n \in B_{\frac{\epsilon}{2}}(w)$. We will now inductively define the rest of the $a_i$'s, two at a time, so that $\hat{X}_{2k+n} \in B_{\frac{\epsilon}{2}}(w)$ for all $k \geq 0$. Suppose we have defined $a_1, \ldots, a_{2k+n}$ properly. Let $a_{2k+n+1} = 0$. Now as $a_{2k+n+2} \to \infty$, $\hat{X}_{2k+n+2} \to \hat{X}_{2k+n}$. Since $\hat{X}_{2k+n} \in B_{\frac{\epsilon}{2}}(w)$, we can choose a sufficiently large $a_{2k+n+2}$ such that $\hat{X}_{2k+n+2} \in B_{\frac{\epsilon}{2}}(w)$ too. The inductive definition is complete.

By the above construction, the sequence $\{a_1, \ldots\}$ contains infinitely many 0's and hence corresponds to a unique point. Thus the largest side length of $\triangle_{a_1 \ldots a_n}$ converges to 0 as $n$ approaches infinity. Hence there exists $N$ such that the longest side length of $\triangle \hat{X}_{2N+n-1} \hat{X}_{2N+n}(X_{2N+n} \hat{+} X_{2N+n-2})$ is less than $\frac{\epsilon}{2}$. But $\hat{X}_{2N+n} \in B_{\frac{\epsilon}{2}}(w)$. Therefore, the partition triangle

$$\triangle \hat{X}_{2N+n-1} \hat{X}_{2N+n}(X_{2N+n} \hat{+} X_{2N+n-2})$$

is contained in $B_\epsilon(w)$. $\square$

Thus every open set contains a partition triangle.

## 7.2 The triangle map is topologically strongly mixing

We will first give the basic definitions and then show that the triangle map is topologically strongly mixing (which implies a number of other dynamical properties). We follow [6] from chapter II, section 4.2.

**Definition 21** *A map* $T : X \to X$ *on a topological space* $X$ *is* topologically strongly mixing *if for any open sets* $U$ *and* $V$ *in* $X$, *there is a positive integer* $N$ *such that for all* $k \geq N$, *we have*

$$T^k U \cap V \neq \emptyset.$$

As discussed in section II.4.4 in [6], topologically strongly mixing implies that the map is topologically ergodic (meaning that given any two open sets $U$ and $V$, there exists some positive integer $N$ such that $T^N U \cap V$ is not empty). The point for us is that topologically strongly mixing is a quite strong condition for topological dynamics.

**Theorem 22** *The triangle map* $T : \triangle \to \triangle$ *is topologically strongly mixing.*

**Proof:**

Recall our notation that a partition triangle $\triangle(a_0, \ldots, a_n)$ denotes all points $(x, y) \in \triangle$ such that

$$
\begin{aligned}
(x, y) &\in \triangle(a_0) \\
T(x, y) &\in \triangle(a_1) \\
T^2(x, y) &\in \triangle(a_2) \\
&\vdots \\
T^n(x, y) &\in \triangle(a_n).
\end{aligned}
$$

Let $U$ and $V$ be two open sets in $\triangle$. By the lemma in the last section, each of these open sets contains a partition triangle.

Denote the partition triangle in $U$ by $\triangle(a_0, \ldots, a_n)$ and the partition triangle in $V$ by $\triangle(b_0, \ldots, b_m)$. Then $\triangle(a_0, \ldots, a_n, b_0, \ldots, b_m)$ is contained in $\triangle(a_0, \ldots, a_n)$ and

$$T^n \triangle (a_0, \ldots, a_n, b_0, \ldots, b_m) = \triangle(b_0, \ldots, b_m).$$

We set $N = n$. For any positive integer $i$, consider $\triangle(a_0, \ldots, a_n, 0, \ldots, 0, b_0, \ldots, b_m)$, where there are $i$ zeros. This partition triangle is contained in $\triangle(a_0, \ldots, a_n)$ and has the property that

$$T^{n+i} \triangle (a_0, \ldots, a_n, 0, \ldots, 0, b_0, \ldots, b_m) = \triangle(b_0, \ldots, b_m).$$

Since $\triangle(a_0, \ldots, a_n, 0, \ldots, 0, b_0, \ldots, b_m)$ is contained in $\triangle(a_0, \ldots, a_n)$, which in turn is contained in the open set $U$ and since $\triangle(b_0, \ldots, b_m)$ is contained in the open set $V$, we must have for all $k \geq N$,

$$T^k U \cap V \neq \emptyset.$$

$\square$

# 8   Future Work

The triangle map has a simple generalization to higher dimensional maps of simplices to themselves, but the corresponding proofs become quite a bit more complicated, at least using the techniques of this paper. Thus one future direction would be to find less cumbersome and more natural arguments for uniqueness and nonuniqueness.

Of course, the main problem is to find an answer to the original Hermite question. For example, is there any way of having the triangle map as a member of a family of maps, each picking up via periodicity a different collection of cubic irrationalities.

Also, the triangle sequence determines a sequence of elements in $\mathbf{SL}(\mathbf{3}, \mathbf{Z})$. It would be interesting to put this in terms of discrete paths in the group $\mathbf{SL}(\mathbf{3}, \mathbf{Z})$ (in particular to relate it to [20].)

## 9 Appendix

This appendix contains derivations of formulae that are used in the earlier parts of this paper. The proofs are straightforward calculations and the formulae themselves give some intuition as to why these ratios and approximations will be useful. Despite this, the precise applicability and usefulness of many of the results of this section can only be seen in context.

### 9.1 Definitions and General Results

The following ratios of side length are important in looking at the behavior of non-unique sequences.

Set

$$\lambda_n = d(\hat{X}_{n-1}, \hat{X}_{n+1})/d(\hat{X}_{n-1}, X_n\hat{+}X_{n-2})$$

$$\lambda'_n = d(X_{n-1}\hat{+}X_{n+1}, \hat{X}_{n+1})/d(\hat{X}_{n-1}, X_n\hat{+}X_{n-2})$$

$$\tilde{\lambda}_n = d(\hat{X}_{n-1}, X_{n-1}\hat{+}X_{n+1})/d(\hat{X}_{n-1}, X_n\hat{+}X_{n-2}).$$

Then

**Proposition 23**

$$\lambda_n = \frac{\frac{x_n + x_{n-2}}{x_{n-1}}}{a_{n+1} + \frac{x_n + x_{n-2}}{x_{n-1}}}$$

$$\tilde{\lambda}_n = \frac{\frac{x_n + x_{n-2}}{x_{n-1}}}{a_{n+1} + 1 + \frac{x_n + x_{n-2}}{x_{n-1}}}$$

$$\lambda_n' = \lambda_n - \tilde{\lambda}_n$$

$$= \frac{\frac{x_n + x_{n-2}}{x_{n-1}}}{(a_{n+1} + \frac{x_n + x_{n-2}}{x_{n-1}})(a_{n+1} + 1 + \frac{x_n + x_{n-2}}{x_{n-1}})}$$

**Proof:** Consider the vector

$$\vec{\mu}_n = (X_n \hat{+} X_{n-2}) - \hat{X}_{n-1}$$

$$= \frac{1}{x_n + x_{n-2}}(X_n + X_{n-2}) - \frac{1}{x_{n-1}}X_{n-1}.$$

Now we have

$$\hat{X}_{n-1}\vec{\hat{X}}_{n+1} = \frac{1}{x_{n+1}}X_{n+1} - \frac{1}{x_{n-1}}X_{n-1}$$

$$= \frac{1}{x_n + x_{n-2} + a_{n+1}x_{n-1}}(X_n + X_{n-2} + a_{n+1}X_{n-1}) -$$

$$\frac{1}{x_{n-1}}X_{n-1}$$

$$= \frac{1}{x_n + x_{n-2} + a_{n+1}x_{n-1}}(X_n + X_{n-2}) -$$

$$\frac{x_n + x_{n-2}}{x_{n-1}(x_n + x_{n-2} + a_{n+1}x_{n-1})}X_{n-1}$$

$$= \frac{x_n + x_{n-2}}{x_n + x_{n-2} + a_{n+1}x_{n-1}}\vec{\mu}_n.$$

Thus

$$\lambda_n = \frac{d(\hat{X}_{n-1}, \hat{X}_{n+1})}{\mu_n} = \frac{x_n + x_{n-2}}{x_n + x_{n-2} + a_{n+1}x_{n-1}},$$

as desired.

Similarly,

$$
\begin{aligned}
\hat{X}_{n-1}(\overrightarrow{X_{n-1}\hat{+}}X_{n+1}) &= \frac{1}{x_{n+1} + x_{n-1}}(X_{n+1} + X_{n-1}) - \frac{1}{x_{n-1}}X_{n-1} \\
&= \frac{1}{x_n + x_{n-2} + (a_{n+1}+1)x_{n-1}}(X_n + X_{n-2} + \\
&\quad (a_{n+1}+1)X_{n-1}) - \frac{1}{x_{n-1}}X_{n-1} \\
&= \frac{1}{x_n + x_{n-2} + (a_{n+1}+1)x_{n-1}}(X_n + X_{n-2}) \\
&\quad - \frac{x_n + x_{n-2}}{x_{n-1}(x_n + x_{n-2} + (a_{n+1}+1)x_{n-1})}X_{n-1} \\
&= \frac{x_n + x_{n-2}}{x_n + x_{n-2} + (a_{n+1}+1)x_{n-1}}\vec{\mu}_n.
\end{aligned}
$$

Then

$$\lambda'_n = \frac{d(X_{n-1}\hat{+}X_{n+1}, \hat{X}_{n+1})}{\mu_n} = \frac{x_n + x_{n-2}}{x_n + x_{n-2} + (a_{n+1}+1)x_{n-1}}.$$

The formula for $\tilde{\lambda}_n$ follows from the above two formulas. $\square$

**Corollary 24** *We have that $\lambda_n$, $\tilde{\lambda}_n$, and $\lambda'_n$ all decrease as $a_{n+1}$ increases.*

## 9.2  Approximations

The following approximations were useful in proving the latter part of the biconditional regarding uniqueness.

**Lemma 25** *Using notation from the previous section, we have*

$$\lambda_n' \leq \frac{1}{\frac{1}{\lambda_n}\left(\frac{1}{\lambda_n} + 1\right)}$$

*and also*

$$\lambda_n' \leq \frac{1}{\frac{1}{\lambda_n}\left(\frac{1}{\lambda_n} - 1\right)}.$$

**Proof:** We have

$$\lambda_n = \frac{\frac{x_n + x_{n-2}}{x_{n-1}}}{a_{n+1} + \frac{x_n + x_{n-2}}{x_{n-1}}}$$

and thus

$$a_{n+1} = \frac{x_n + x_{n-2}}{x_{n-1}}\left(\frac{1}{\lambda_n} - 1\right) \geq \frac{1}{\lambda_n} - 1.$$

Hence

$$
\begin{aligned}
\lambda_n' &= \frac{\frac{x_n + x_{n-2}}{x_{n-1}}}{(a_{n+1} + \frac{x_n + x_{n-2}}{x_{n-1}})(a_{n+1} + 1 + \frac{x_n + x_{n-2}}{x_{n-1}})} \\
&= \frac{\lambda_n}{a_{n+1} + 1 + \frac{x_n + x_{n-2}}{x_{n-1}}} \\
&\leq \frac{\lambda_n}{\left(\frac{1}{\lambda_n} - 1\right) + 1 + 1} \\
&= \frac{1}{\frac{1}{\lambda_n}\left(\frac{1}{\lambda_n} + 1\right)}.
\end{aligned}
$$

Similarly,

$$\tilde{\lambda}_n = \frac{\frac{x_n + x_{n-2}}{x_{n-1}}}{a_{n+1} + 1 + \frac{x_n + x_{n-2}}{x_{n-1}}}$$

implies that

$$a_{n+1} = \frac{x_n + x_{n-2}}{x_{n-1}}\left(\frac{1}{\tilde{\lambda}_n} - 1\right) - 1$$

$$\geq 1\left(\frac{1}{\lambda_n} - 1\right) - 1$$
$$= \frac{1}{\lambda_n} - 2.$$

Hence a similar argument as above yields the second desired inequality. □

**Lemma 26** *We have*

$$\tau_{n+1} \geq \tau_n - \lambda_n(\rho_n + \tau_n)$$

*and*

$$\rho_{n+1} \leq \frac{\rho_n + \tau_n}{\frac{1}{\lambda_n}\left(\frac{1}{\lambda_n} + 1\right)}.$$

**Proof:** The length of the vector $\vec{\tau}_{n+1}$ is:

$$\begin{aligned}
\tau_{n+1} &= |\lambda_n\vec{\rho}_n - (1 - \lambda_n)\vec{\tau}_n| \\
&\geq |\lambda_n\vec{\rho}_n| - |(1 - \lambda_n)\vec{\tau}_n| \\
&= \lambda_n\rho_n - (1 - \lambda_n)\tau_n \\
&= \tau_n - \lambda_n(\rho_n + \tau_n).
\end{aligned}$$

Second,

$$\begin{aligned}
\rho_{n+1} &= \lambda_n' d(X_n \hat{+} X_{n-2}, \hat{X}_{n-1}) \\
&\leq \lambda_n'\left(d(X_n \hat{+} X_{n-2}, \hat{X}_n) + d(\hat{X}_n, \hat{X}_{n-1})\right) \\
&= \lambda_n'(\rho_n + \tau_n) \\
&\leq \frac{\rho_n + \tau_n}{\frac{1}{\lambda_n}\left(\frac{1}{\lambda_n} + 1\right)}.
\end{aligned}$$

□